\magnification=1000
\hsize=11.7cm
\vsize=18.9cm
\lineskip2pt \lineskiplimit2pt
\nopagenumbers

\hoffset=-1truein
\voffset=-1truein

\advance\voffset by 4truecm
\advance\hoffset by 4.5truecm

\newif\ifentete

\headline{\ifentete\ifodd	\count0 
      \rlap{\head}\hfill\tenrm\llap{\the\count0}\relax
    \else
        \tenrm\rlap{\the\count0}\hfill\llap{\head} \relax
    \fi\else
\global\entetetrue\fi}

\def\entete#1{\entetefalse\gdef\head{#1}}
\entete{}

\input amssym.def
\input amssym.tex

\def\-{\hbox{-}}
\def\.{{\cdot}}
\def\O{{\cal O}}

\def\int{\frak i\frak n\frak t}

\def\qq{\quad{\rm and}\quad}

\def\too{\longrightarrow}

 3
 2
\font\large=cmr10  scaled \magstep 2
 2
 2
 2
\font\cds=cmr7

\centerline{\large Nilpotent extensions of blocks}

\vskip 0.5cm

\centerline{\bf Lluis Puig }
\medskip
\centerline{\cds CNRS, Institut de Math\'ematiques de Jussieu}
\smallskip
\centerline{\cds 6 Av Bizet, 94340 Joinville-le-Pont, France}
\smallskip
\centerline{\cds puig@math.jussieu.fr}

\vskip 0.5cm
\noindent
{\bf £1\phantom{.} Introduction}
\bigskip

£1.1\phantom{.} The {\it nilpotent blocks\/} over an algebraically closed
field of characteristic $p > 0$ were introduced in [2] as a
translation for blocks of the well--known Frobenius Criterion
on $p$-nilpotency for finite groups. They correspond to
the simplest situation with respect to the so--called {\it
fusion\/} inside a defect group, and the structure of the source algebras
of the nilpotent blocks determined in~[9] confirms that these blocks represent indeed the easiest
possible situation.

\medskip
£1.2\phantom{.} However, when the field of coefficients is not algebraically closed, together with Fan Yun
we have seen in [3] that, in the general situation, the structure of the source algebra of a block which, 
after a suitable scalar extension, decomposes in a sum of nilpotent blocks --- a structure that we determine in [3]
--- need not be so simple
.

 \medskip
£1.3\phantom{.} At that time, we already knew some examples of a similar fact  in group extensions,
namely that a {\it non-nilpotent\/} block of a normal subgroup $H$ of a finite group $G$ may decompose 
in a sum of nilpotent blocks of $G\,.$ In this case, we also have been able to describe the 
source algebra structure, which is quite similar to (but easier than)  the structure described in [3]. With a big delay,
we explain this result here.

\medskip
£1.4\phantom{.} Actually, this phenomenon  is perhaps better described by saying that {\it a normal sub-block
of a nilpotent block need not be nilpotent\/}. However, the normal sub-blocks of nilpotent blocks are quite special:
 they are {\it basically Morita equivalent\/} [15,~\S7] to the corresponding block of their {\it inertial subgroup\/}.
Then, as a matter of fact, a normal sub-block of such a block still fulfills the same condition.

\medskip
£1.5\phantom{.} Thus, let us call {\it inertial block\/} any block of a finite group that is {\it basically
Morita equivalent\/}  [15,~\S7] to the corresponding block of its {\it inertial subgroup\/}; as a matter of fact,
in [12, Corollaire~3.6] we already exhibit a large family of inertial blocks; see also [14, Appendix]. The main purpose of this paper
is to prove that {\it a normal sub-block of an inertial block is again an inertial block\/}. Since a nilpotent block
is {\it basically Morita equivalent\/} to its defect group [9, Theorem~1.6 and (1.8.1)], and the corresponding block of its {\it inertial subgroup\/} is also nilpotent, a nilpotent block is, in particular, an inertial block and thus, our main result applies.
\eject

\bigskip
\bigskip
\noindent
{\bf £2\phantom{.} Quoted results and inertial blocks}
\bigskip

£2.1\phantom{.} Throughout this paper $p$ is a fixed prime number, $k$ an algebraically closed field of 
characteristic $p$ and $\O$ a complete discrete valuation ring of characteristic zero having the {\it residue field\/} $k\,.$
Let $G$ be a finite group; following Green [5], a {\it $G\-$algebra\/} is a torsion-free $\O\-$algebra $A$ of finite 
$\O\-$rank endowed with a $G\-$action; we say that $A$ is {\it primitive\/} if the unity element is
primitive in~$A^G\,.$ A  $G\-$algebra homomorphism from~$A$ to another   $G\-$algebra~$A'$ is a 
{\it not necessarily unitary\/} algebra homomorphism $f\,\colon A\to A'$ compatible
with the $G\-$actions. We say that $f$ is an {\it embedding\/} whenever  
$${\rm Ker}(f) = \{0\}\qq {\rm Im}(f) = f(1_A)A'f(1_A)
\eqno £2.1.1,$$
and that $f$ is a {\it strict semicovering\/} if $f$ is {\it unitary\/},  the {\it radical\/} $J(A)$ of $A$ contains
 ${\rm Ker}(f)$ and, for any $p\-$subgroup $P$ of $G\,,$ $J(A'^P)$ contains $f\big(J(A^P)\big)$  and
$f(i)$ is primitive in~$A'^P$ for any primitive idempotent $i$ of $A^P$~[6,~\S3].

\medskip
£2.2\phantom{.} Recall that, for any subgroup $H$ of $G\,,$
a {\it point\/} $\alpha$ of $H$ on $A$ is an $(A^H)^*\-$conjugacy class
of primitive idempotents of $A^H$ and the pair $H_\alpha$ is a {\it
pointed group\/} on $A$ [7, 1.1]; if $H = \{1\}\,,$ we simply say that $\alpha$ is a {\it point\/} of~$A\,.$ For any $i\in \alpha\,,$ $iAi$ has an evident
structure of  $H\-$algebra and we denote by~$A_\alpha$ one of these mutually
$(A^H)^*\-$conjugate $H\-$algebras and by~$A(H_\alpha)$ the {\it simple
quotient\/} of $A^H$ determined by $\alpha\,;$ we call {\it multiplicity\/} of $\alpha$ the {\it square root\/} of the
dimension of $A(H_\alpha)\,.$ If  $f\,\colon A\to A'$ is a $G\-$algebra homomorphism
and $\alpha'$ a point of $H$ on $A'\,,$ we call {\it multiplicity\/} ${\rm m}(f)_\alpha^{\alpha'}$ of $f$ at 
$(\alpha,\alpha')$ the dimension of the image of $f(i)A'^Hi'$ in $A'(H_{\alpha'})$ for $i\in \alpha$ and  
$i'\in \alpha'\,;$ we still consider the $H\-$algebra  $A'_\alpha = f(i)A'f(i)$
together with the unitary $H\-$algebra homomorphism induced by $f$ and the embedding of $H\-$algebras
$$A_\alpha\too A'_\alpha \longleftarrow A'_{\alpha'}
\eqno £2.2.1.$$
A second pointed group $K_\beta$ on $A$ is {\it contained\/} in~$H_\alpha$ if $K\i H$ and, 
for any $i\in\alpha\,,$ there is $j\in \beta$ such that~[7,~1.1]
$$ij = j = ji
\eqno £2.2.2;$$
then, it is clear that the $(A^K)^*\-$conjugation induces $K\-$algebra
embeddings
$$f_\beta^\alpha : A_\beta\too {\rm Res}^{H}_{K} (A_\alpha)
\eqno £2.2.3.$$

\medskip
£2.3\phantom{.} Following Brou\'e, for any $p\-$subgroup $P$ of $G$ we consider the {\it Brauer quotient\/}
and the {\it Brauer homomorphism\/} [1, 1.2]
$${\rm Br}^A_P : A^P\too A (P) =  A^P\Big/\sum_Q A^P_Q
\eqno £2.3.1,$$
where $Q$ runs over the set of proper subgroups of $P\,,$ and call  {\it local\/} any point~$\gamma$ of $P$ 
on $A$ not contained in   ${\rm Ker(Br}^A_P)$ [7, 1.1]. Recall that
{\it a local pointed group $P_\gamma$ contained in $H_\alpha$ is maximal if and
only if ${\rm Br}_P(\alpha)\i A(P_\gamma)^{N_H (P_\gamma)}_P$\/}\break
\eject
\noindent
[7, Proposition~1.3] and then {\it the
$P\-$algebra~$A_\gamma$\/} --- called a {\it source algebra\/} of $A_\alpha$ --- {\it
is Morita equivalent to~$A_\alpha$\/} [17,~6.10]; moreover,
{\it the maximal local pointed groups~$P_\gamma$ contained in $H_\alpha$\/} ---
called the {\it defect pointed groups\/} of $H_\alpha$ --- {\it are  mutually
$H\-$conjugate\/} [7,~Theorem~1.2].

\medskip
£2.4\phantom{.} Let us say that $A$ is a {\it $p\-$permutation
$G\-$algebra\/} if a Sylow $p\-$subgroup of $G$ stabilizes a basis of $A$
[1, 1.1]. In this case, recall that if $P$ is a $p\-$subgroup of $G$ and $Q$ a normal subgroup of $P$
then the corresponding Brauer homomorphisms induce a $k\-$algebra isomorphism [1, Proposition~1.5]
$$\big(A(Q)\big)(P/Q)\cong A(P)
\eqno £2.4.1;$$
moreover, choosing a point $\alpha$ of $G$ on~$A\,,$ we call 
{\it Brauer $(\alpha,G)\-$pair\/} any pair $(P,e_A)$ formed by a
$p\-$subgroup~$P$ of $G$ such that ${\rm Br}^A_P (\alpha)\not= \{0\}$ and by
 a primitive idempotent $e_A$ of the center $Z\big( A (P)\big)$ of $A(P)$ such
that 
$$e_A\. {\rm Br}^A_P (\alpha)\not= \{0\}
\eqno £2.4.2;$$
 note that any local pointed group $Q_\delta$ on $A$ {\it contained\/} in $G_\alpha$
determines a Brauer $(\alpha,G)\-$pair $(Q,f_A)$ fulfilling $f_A \.{\rm Br}^A_Q
(\delta)\not= \{0\}\,.$

\medskip
£2.5\phantom{.} Then, it follows from Theorem~1.8 in [1] that {\it the
inclusion between the local pointed groups on $A$  induces an inclusion
between the Brauer $(\alpha,G)\-$pairs\/}; explicitly, if $(P,e_A)$ and
$(Q,f_A)$ are two Brauer $(\alpha,G)\-$pairs then we have
$$(Q,f_A)\i (P,e_A)
\eqno £2.5.1\phantom{.}$$
whenever there are local pointed groups $P_\gamma$ and $Q_\delta$ on $A$
fulfilling $$Q_\delta\i P_\gamma\i G_\alpha\quad ,\quad f_A \.{\rm Br}^A_Q (\delta)\not=
\{0\} \qq    e_A\. {\rm Br}^A_P (\gamma)\not= \{0\}
\eqno £2.5.2.$$
Actually, according to the same result, for any $p\-$subgroup $P$ of $G\,,$
any primitive idempotent $e_A$ of $Z\big(A(P)\big)$ fulfilling $e_A\. {\rm Br}^A_P (\alpha)\not= \{0\}$ 
and any subgroup~$Q$ of~$P\,,$ there is a unique
primitive idempotent $f_A$ of $Z\big(A(Q)\big)$ fulfilling 
$$e_A\. {\rm Br}^A_P (\alpha)\not=\{0\}\qq (Q,f_A)\i (P,e_A)
\eqno £2.5.3.$$
Once again, {\it the maximal Brauer $(\alpha,G)\-$pairs are pairwise
$G\-$conjugate\/} [1, Theorem~1.14].

\medskip
£2.6\phantom{.} Here, we are specially interested in the $G\-$algebras $A$ endowed with a group
homomorphism $\rho\,\colon G\to A^*$ inducing the action of $G$ on $A\,,$ called {\it $G\-$interior algebras\/};
 in this case, for any pointed group $H_\alpha$ on $A\,,$ $A_\alpha = iAi$ has a  structure of 
{\it $H\-$interior algebra\/} mapping $y\in H$ on $\rho (y)i = i\rho(y)\,;$ moreover, setting $x\.a\.y = \rho(x)a\rho(y)$
for any $a\in A$ and any $x,y\in G\,,$ a  $G\-$interior algebra homomorphism 
from~$A$ to another   $G\-$interior algebra~$A'$ is a $G\-$algebra homomorphism $f\,\colon A\to A'$ fulfilling
$$f(x\.a\.y) = x\.f(a)\.y
\eqno £2.6.1.$$
\eject

\medskip
£2.7\phantom{.} In particular, if $H_\alpha$ and $K_\beta$ are two pointed groups on $A\,,$ we say that an
injective group homomorphism $\varphi\,\colon K\to H$ is an
{\it $A\-$fusion from $ K_\beta$ to~$H_\alpha$\/} whenever there is
a  $K\-$interior algebra~{\it embedding\/} 
$$f_{\varphi} : A_\beta\too {\rm Res}^{H}_{K} (A_\alpha) 
\eqno £2.7.1\phantom{.}$$
such that the inclusion $A_\beta\i A$ and the composition of $f_{\varphi}$
with the inclusion $A_\alpha\i A$ are $A^*\-$conjugate; we denote by $F_A
( K_\beta,H_\alpha)$ the set of $H\-$conjugacy classes of  $A\-$fusions from $ K_\beta$ to~$H_\alpha$
and, as usual, we write $F_A (H_\alpha)$ instead of $F_A
(H_\alpha,H_\alpha)\,.$ If $A_\alpha = iAi$ for $i\in \alpha\,,$  it follows
from [8, Corollary~2.13] that we have a group homomorphism
$$F_A (H_\alpha)\too N_{A_\alpha^{^*}}(H\.i)\big/H\.(A_\alpha^H)^*
\eqno £2.7.2\phantom{.}$$
and if $H$ is a $p\-$group then we consider the $k^*\-$group $\hat F_A (H_\alpha)$
defined by the {\it pull-back\/}
$$\matrix{F_A (H_\alpha)&\too& N_{A_\alpha^{^*}}(H\.i)/H\.(A_\alpha^H)^*\cr 
\uparrow&\phantom{\big\uparrow}&\uparrow\cr
\hat F_A (H_\alpha) &\too &N_{A_\alpha^{^*}}(H\.i)\big/H\.\big(i + J(A_\alpha^H)\big)\cr}
\eqno £2.7.3.$$

\medskip
£2.8\phantom{.} Recall that, for any subgroup $H$ of $G$ and any $H\-$interior algebra~$B\,,$ the {\it induced $G\-$interior
algebra\/} is the induced bimodule
$${\rm Ind}_{ H}^{ G}(B) = k G\otimes_{k H} B
\otimes_{k  H} k  G
\eqno £2.8.1,$$
endowed with the distributive product defined by the {\it formula\/}
$$( x\otimes b\otimes  y)( x'\otimes b'\otimes  y')
=\cases{ x\otimes b. y x'.b'\otimes  y'&if $ y x'\in
 H$\cr  {}&{}\cr 
0 &otherwise\cr}
\eqno £2.8.2\phantom{.}$$
where $ x, y, x', y'\in  G$ and $b,b'\in B\,,$ and with
the structural homomorphism
$$ G\too {\rm Ind}_{ H}^{ G}(B)
\eqno £2.8.3\phantom{.}$$
mapping $ x\in  G$ on the element
$$\sum_{ y} x y\otimes 1_B\otimes  y^{-1} = \sum_{ y} y\otimes 1_B
\otimes  y^{-1} x
\eqno £2.8.4\phantom{.}$$
where $ y\in  G$ runs over a set of representatives for 
$ G/ H\,.$

\medskip
£2.9\phantom{.} Obviously, the {\it group algebra\/} $\O G$ is a
$p\-$permutation $G\-$interior algebra and, for any primitive idempotent $b$ of
$Z(\O G)$ --- called an {\it $\O\-$block\/} of~$G$ -- the conjugacy class $\alpha = \{b\}$ is a {\it point\/}
 of $G$ on $\O G\,.$ Moreover, for any $p\-$subgroup $P$ of $G\,,$ the Brauer homomorphism 
 ${\rm Br}_P = {\rm Br}_P^{kG}$ induces a $k\-$algebra isomorphism [10, 2.8.4]
$$kC_G (P)\cong (\O G)(P)
\eqno £2.9.1;$$
\eject
\noindent
thus, up to identification throughout this isomorphism, in a Brauer $(\{b\},G)\-$ pair $(P,e)$  as defined above
--- called {\it Brauer $(b,G)\-$pair\/} from now on --- $e$ is nothing but a $k\-$block of $C_G (P)$  such that $e{\rm Br}_P (b)\not= 0\,.$ Setting 
$$\bar C_G (P) = C_G (P)/Z(P)
\eqno £2.9.2,$$
recall that the image  $\bar e$ of $e$ in $k \bar C_G (P)$ is a $k\-$block of $\bar C_G (P)$ and that 
the {\it Brauer First Main Theorem\/} affirms that {\it $(P,e)$ is maximal if and only if the $k\-$algebra 
$k\bar C_G (P)\bar e$  is simple and the inertial
quotient 
$$E = N_G (P,e)/P\.C_G (P)
\eqno £2.9.3\phantom{.}$$
 is a $p'\-$group\/}~[17,~Theorem~10.14].

\medskip
£2.10\phantom{.} For any
$p\-$subgroup $P$ of~$G$ and any subgroup $H$ of $N_G (P)$ containing $P\.C_G(P)\,,$
we have
$${\rm Br}_P\big((\O G)^H\big) = (\O G)(P)^H
\eqno £2.10.1\phantom{.}$$
and therefore {\it any $k\-$block $e$ of $C_G (P)$ determines a unique point~$\beta$ of $H$ on~$\O G$ 
{\rm (cf.~£2.2)} such that $H_\beta$ contains $P_\gamma$ for a local point $\gamma$ of $P$ on $\O G$
fulfilling\/}~[9,~Lemma~3.9]
$$e\. {\rm Br}_P (\gamma)\not= \{0\}
\eqno £2.10.2.$$
Recall that, if $Q$ is a subgroup of $P$ such that $C_G (Q)\i H$ then the
$k\-$blocks of $C_G (Q) = C_H(Q)$ determined by $(P,e)$ from $G$ and from $H$ coincide
[1,~Theorem~1.8].  Note that if $P$ is normal in $G$ then the kernel of the
obvious $k\-$algebra homomorphism $kG\to k(G/P)$ is contained in the {\it
radical\/} $J(kG)$ and contains~${\rm Ker}({\rm Br}_P)\,;$
thus, in this case, isomorphism~£2.9.1 implies that  {\it any point of~$P$ on $kG$ is local.\/}

\medskip
£2.11\phantom{.}  Moreover, for any local pointed group $P_\gamma$ on $\O G\,,$ the action of $N_G (P_\gamma)$
on the simple algebra $(\O G)(P_\gamma)$ (cf.~£2.2) determines  a central $k^*\-$ex-tension or, equivalently, 
a {\it $k^*\-$group\/} $\hat N_G (P_\gamma)$ [10, \S5] and it is clear that the Brauer homomorphism ${\rm Br}_P$ determines a $N_G(P_\gamma)\-$stable  injective group homomorphism from $C_G (P)$ to $\hat N_G (P_\gamma)\,.$ Then, up to a suitable identification, we set  
$$E_G (P_\gamma) = N_G (P_\gamma)/P\.C_G (P)\qq \hat E_G (P_\gamma) = \hat N_G (P_\gamma)/P\.C_G (P)
  \eqno £2.11.1;$$
recall that from [8, Theorem~3.1] and [10,~Proposition~6.12] we obtain a {\it canonical\/} $k^*\-$group 
isomorphism (cf.~£2.7.3)
$$\hat E_G (P_\gamma)^\circ\cong \hat F_{\O G}(P_\gamma)
\eqno £2.11.2.$$

\medskip
£2.12\phantom{.} In particular, a maximal local pointed group $P_\gamma$ on $\O Gb$ determines a
 $k\-$block $e$ of $C_G (P)\,,$ which is still a $k\-$block of the group 
 $$N = N_G (P_\gamma) = N_G (P,e)
 \eqno £2.12.1\phantom{.}$$
--- called the {\it inertial subgroup\/} of $b$ --- and also determines  a unique 
 point $\nu$ of~$N$ on $\O Gb$ such that $P_\gamma\i N_\nu$ (cf.~£2.10); obviously, we have 
 $E = E_G (P_\gamma)$ (cf.~£2.9.3), $P_\gamma$ is still
 a {\it defect pointed group\/} of $N_\nu$ and $(P,e)$ is a maximal Brauer $(\hat e,N)\-$pair, where $\hat e$
 denotes the  $\O\-$block of~$N$ lifting $e\,.$ As above, $N$ acts on the simple $k\-$algebra (cf.~£2.9)
 $$k\bar C_G (P)\bar e\cong (\O G)(P_\gamma)
 \eqno £2.12.2\phantom{.}$$
  and therefore we get {\it $k^*\-$groups\/} $\hat N$ and $\hat E^\circ = \hat E_G (P_\gamma)\,.$

 \medskip
£2.13\phantom{.} Moreover, since $E$ is a $p'\-$group, it follows from [17, Lemma~14.10] that the short exact sequence
$$1\too P/Z(P) \too N/C_G (P)\too E\too 1
  \eqno £2.13.2\phantom{.}$$  
splits and that all the splitings are conjugate to each other; thus, any spliting determines an action of   $E$
 on $P$ and it is easily checked that the semidirect products
 $$L = P\rtimes E\qq \hat L = P\rtimes \hat E
 \eqno £2.13.3\phantom{.}$$
do not depend on our choice. At this point,  it follows from  [10, Proposition~14.6] that the source algebra 
  of the block $\hat e$ of $N$ is isomorphic to the  $P\-$interior algebra $\O_*\hat L\,,$ and therefore it follows from [3, Proposition~4.10] that the multiplication in $\O Gb$ by a suitable idempotent $\ell\in \nu$ determines an injective unitary $P\-$interior algebra homomorphism
$$\O_*\hat L\too (\O G)_\gamma
\eqno £2.13.4.$$

\medskip
£2.14\phantom{.} On the other hand, a {\it Dade $P\-$algebra\/} over $\O$
is a $p\-$permutation $P\-$algebra $S$ which is a {\it full matrix algebra over $\O$\/} and
fulfills $S(P)\not=\{0\}$ [11, 1.3].  For any subgroup $Q$ of~$P\,,$ setting $\bar N_P (Q) = N_P (Q)/Q$ 
we have (cf.~£2.4.1)
$$\big(S(Q)\big)\big(\bar N_P(Q)\big)\cong S\big(N_P(Q)\big)
\eqno £2.14.1\phantom{.}$$
 and therefore ${\rm Res}^P_Q (S)$ is
a Dade $Q\-$algebra; moreover, it follows from [11,~1.8] that the {\it Brauer
quotient\/} $S(Q)$ is a Dade $\bar N_P (Q)\-$algebra over $k\,;$ thus, $Q$ has a unique {\it local point\/} on $S\,.$
In particular, if $S$ is {\it primitive\/} (cf.~£2.1) then $S(P)\cong k$ and therefore we have
$${\rm dim}(S)\equiv 1\pmod{p}\eqno £2.14.2,$$
so that the action of $P$ on $S$ can be lifted to a unique group homomorphism from $P$ to the kernel
of the determinant ${\rm det}_S$ over $S\,;$ at this point, it follows from [11, 3.13] that  the action of $P$ on $S$
always can be lifted to a well-determined $P\-$interior algebra structure for $S\,.$
\eject

\medskip
£2.15\phantom{.}  Recall that a block $b$ of $G$ is called {\it nilpotent\/} whenever the quotients 
$N_G (Q,f)/C_G (Q)$ are $p\-$groups for all the Brauer $(b,G)\-$pairs $(Q,f)$ [2, De-finition~1.1]; by the
main result in~[9], {\it the block $b$ is nilpotent if and only if, for a
maximal local pointed group $P_\gamma$ on $\O Gb\,,$ $P$ stabilizes a
unitary primitive Dade  $P\-$subalgebra $S$ of $(\O Gb)_\gamma$ fulfilling\/}
$$(\O Gb)_\gamma = SP\cong S\otimes_\O \O P
\eqno £2.15.1\phantom{.}$$ 
where we denote by $SP$ the obvious $\O\-$algebra $\bigoplus_{u\in P}Su$ and,
for the right-hand isomorphism, we consider the  well-determined $P\-$interior algebra structure for $S\,.$

\medskip
£2.16\phantom{.} Now, with the notation in~£2.12 above, we say that the block $b$ of~$G$ is {\it inertial\/} 
if it is {\it basically Morita equivalent\/} [15, 7.3] to the corresponding block~$\hat e$ of the
{\it  inertial subgroup\/} $N$ of~$b$ or, equivalently, if there is a primitive Dade $P\-$algebra $S$ such that we have
a $P\-$interior algebra embedding [15, Theorem~6.9 and Corollary~7.4]
$$(\O G)_\gamma\too S\otimes_\O \O_*\hat L
\eqno £2.16.1.$$
Note that, in this case, in fact {\it we  have a $P\-$interior algebra isomorphism
$$(\O G)_\gamma\cong S\otimes_\O \O_*\hat L
\eqno £2.16.2\phantom{.}$$
and the Dade $P\-$algebra $S$ is uniquely determined\/}; indeed, the 
uniqueness of~$S$ follows from [19, Lemma~4.5] and it is easily checked that 
$$(S\otimes_\O \O_*\hat L)(P)\cong S(P)\otimes_k (\O_*\hat L)(P)\cong kZ(P)
\eqno £2.16.3\phantom{.}$$
and that the kernel of the Brauer homomorphism ${\rm Br}^{S\otimes_\O \O_*\hat L}_P$ is contained
in the radical of~$S\otimes_\O \O_*\hat L\,,$ so that this  $P\-$interior algebra is also primitive.

\bigskip
\bigskip
\noindent
{\bf £3\phantom{.} Normal sub-blocks of inertial blocks}
\bigskip

£3.1\phantom{.} Let~$G$ be a  finite group, $b$ an {\it $\O\-$block\/} of $G$ and
$(P,e)$ a maximal Brauer $(b,G)\-$pair (cf.~£2.9). Let us say that an {\it $\O\-$block\/} $c$ of 
a normal subgroup~$H$ of $G$ is a {\it normal sub-block\/} of $b$ if  we have $cb\not= 0\,;$
we are interested in the relationship between the source algebras of $b$ and $c\,,$ specially in the case
where $b$ is {\it inertial\/}.

\medskip
£3.2\phantom{.}  Note that  we have $b{\rm Tr}_{G_c}^{G}(c) = b$ where $G_c$ denotes the
stabilizer of $c$ in $G\,;$ since we know that $e{\rm Br}_P (b)\not= 0$ (cf.~£2.9), up to modifying our
choice of $(P,e)$ we may assume that $P$ stabilizes $c\,;$ then, considering the $G\-$stable
semisimple $k\-$subalgebra $\sum_x k\.bc^x$ of $k G\,,$ where $x\in G$ runs over a set of representatives for
$G/G_c\,,$ it follows from [19, Proposition~3.5] that $bc$ is an $\O\-$block of $G_c$ and
 that $P$ remains a defect $p\-$subgroup of this block, and then from [19, Proposition~3.2] that we have
$$\O Gb\cong {\rm Ind}_{G_c}^G (\O G_c\, bc)
\eqno £3.2.1,$$
so that the source algebras of the $\O\-$block $b$ of $G$ and of the block $bc$ of $G_c$ are
isomorphic.

\medskip
£3.3\phantom{.} Thus, from now on we assume that $G$ fixes~$c\,,$ so that we have $bc = b\,.$ 
Then, note that $\alpha = \{c\}$ is a point of $G$ on $\O H$ (cf.~£2.2), so that, choosing
 a block $e^{_H}$ of~$C_{H} (P)$ such that $e^{_H}e\not= 0\,,$ $(P,e^{_H})$ is a {\it Brauer $(\alpha,G)\-$pair\/}
(cf.~£2.4 and~£2.9.1) and it follows from the proof of [18, Proposition~15.9] that we may choose 
a maximal Brauer $(c,H)\-$pair $(Q,f^{^H})$ fulfilling
$$(Q,f^{^H})\i (P,e^{_H})\quad ,\quad Q = H\cap P\qq e{\rm Br}_P(f^{^H})\not= 0
\eqno £3.3.1.$$
Now, denote by $\gamma^{_G}$ and $\delta$ the respective local points of $P$ and $Q$ on $\O G$ and $\O H$
determined by $e$ and $f^{^H}\,;$ as above, let us denote by $F$ the {\it inertial quotient\/} of $c\,;$ that is to say, we set (cf.~£2.9 and~£2.11)
$$F  =  E_H (Q_\delta) = F_{\O H}(Q_\delta)\qq \hat F =\hat  E_H (Q_\delta)^\circ\cong \hat F_{\O H}(Q_\delta)\eqno £3.3.2.$$

\medskip
£3.4\phantom{.} Since we have $e{\rm Br}_P(f^{^H})\not= 0$ and $f^{^H}$ is
$P\-$stable, from the obvious commutative diagram
$$\matrix{(\O H)(Q)&\too &(\O G)(Q)\cr
\cup&&\cup\cr
(\O H)(Q)^P&\too & (\O G)(Q)^P\cr
\downarrow&\phantom{\Big\downarrow}&\downarrow\cr
(\O H)(P)&\too &(\O G)(P)\cr}
\eqno £3.4.1\phantom{.}$$
we get a local point $\delta^{^G}$ of $Q$ on $\O G$ such that the multiplicity ${\rm m}_\delta^{\delta^{^G}}\!$
of the inclusion $(\O H)^Q\i (\O G)^Q$ at $(\delta,\delta^{^G})$ (cf.~£2.2) is not zero and $Q_{\delta^{^G}}$ is contained  in~$P_{\gamma^{^G}}\,;$ similarly, we get a local point $\gamma$ of $P$ on $\O H$ fulfilling
$${\rm m}_\gamma^{\gamma^{_G}}\not= 0\qq Q_\delta\i P_\gamma
\eqno £3.4.2\phantom{.}$$
At this point,  the following commutative diagram (cf.~2.2.1)
$$\matrix{&&\hskip-40pt{\rm Res}_Q^P(\O H)_\gamma & \hskip-20pt\too &  \hskip-10pt{\rm Res}_Q^P(\O G)_\gamma\cr
&\hskip-30pt\nearrow&&\hskip-10pt\nearrow&\hskip-10pt\uparrow\cr
(\O H)_\delta&\too & (\O G)_\delta&&\hskip-10pt{\rm Res}_Q^P(\O G)_{\gamma^{^G}}\cr
&&\uparrow&\hskip-10pt\nearrow\cr
&&(\O G)_{\delta^{^G}}\cr}
\eqno £3.4.3,$$
where all the $Q\-$interior algebra homomorphisms but the horizontal ones are embeddings,
 already provides some relationship between the source algebras of 
$b$ and $c$ (cf.~£2.2). 
\eject

\medskip
£3.5\phantom{.} If $R_\varepsilon$ is a local pointed group on $\O H\,,$ we set 
$$C_G (R_\varepsilon) = C_G (R)\cap N_G (R_\varepsilon)\qq  E_G (R_\varepsilon) = 
N_G (R_\varepsilon)/R\.C_G (R_\varepsilon)
\eqno £3.5.1\phantom{.}$$
and denote by $b(\varepsilon)$ the block of $C_H (R)$ determined by $\varepsilon\,,$ and by 
$\bar b(\varepsilon)$ the image of $b(\varepsilon)$ in $k\bar C_H (R) = k\big(C_H (R)/Z(R)\big)\,;$ recall that we have 
a canonical
$\bar C_G (R)\-$interior algebra isomorphism [19, Proposition~3.2]
$$k\bar C_G (R){\rm Tr}_{\bar C_G (R_\varepsilon)}^{\bar C_G (R)} \big(\bar b (\varepsilon)\big) \cong
{\rm Ind}_{\bar C_G (R_\varepsilon)}^{\bar C_G (R)} \big(k\bar C_G (R_\varepsilon)\bar b(\varepsilon)\big)
\eqno £3.5.2.$$
Moreover, note that if $\varepsilon^{_G}$ is a local point of $R$ on $\O G$ such that 
${\rm m}_\varepsilon^{\varepsilon^{_G}}\not=\! 0$ then we have
$$E_G (R_{\varepsilon^{_G}})\i E_G (R_\varepsilon)
\eqno £3.5.3;$$
indeed, the restriction to $C_H (R)$ of a simple $kC_G (R)\-$module determined by~$\varepsilon^{_G}$
is semisimple (cf.~£2.9.1) and therefore $C_G (R)$ acts transitively on the set of local points $\varepsilon'$ 
of $R$ on $\O H$ such that ${\rm m}_{\varepsilon'}^{\varepsilon^{_G}}\not=\! 0\,,$ so that we have
$$N_G (R_{\varepsilon^{_G}}) \i C_G(R)\.N_G (R_\varepsilon)
\eqno £3.5.4.$$
Then, we also define $E_H(R_{\varepsilon^{_G}}) = N_H (R_{\varepsilon^{_G}})/R\.C_H(R)\,.$

\medskip
£3.6\phantom{.} Since $(Q,f^{^H})$ is a maximal Brauer $(c,H)\-$pair, we have (cf.~£2.12.2)
$$k\bar C_H (Q)\bar f^{^H}\cong (\O H)(Q_\delta)
\eqno £3.6.1\phantom{.}$$
and, according to the very definition of the $k^*\-$group $\hat N_G (Q_\delta)\,,$ we also have a $k^*\-$group homomorphism
$$\hat N_G (Q_\delta)\too \big(k\bar C_H (Q)\bar f^{^H}\big)^*
\eqno £3.6.2;$$
then, denoting by $\hat C_G (Q_\delta)$ the corresponding $k^*\-$subgroup of 
$\hat N_G (Q_\delta)$ and setting
$$Z =  C_G (Q_\delta)/C_H (Q)\qq \hat Z = \hat  C_G (Q_\delta)/C_H (Q)
\eqno £3.6.3,$$
 it follows from  [19, Theorem~3.7] that we have a canonical 
$\bar C_G (Q_\delta)\-$interior algebra isomorphism
$$k\bar C_G (Q_\delta)\bar f^{^H}\!\cong k\bar C_H (Q)\bar f^{^H}\otimes_k 
(k_*\hat Z)^\circ
\eqno £3.6.4.$$
Now, this isomorphism and the corresponding isomorphism~£3.5.2 determine a $k\-$algebra isomorphism
$$Z\big(k\bar C_G (Q)\big){\rm Tr}_{\bar C_G (Q_\delta)}^{\bar C_G (Q)} (\bar f^{^H})
\cong Z (k_*\hat Z)
\eqno £3.6.5,$$
and induce a bijection between the set of local points   $\delta^{_G}$ of $Q$ on $\O Gb$ such that ${\rm m}_{\delta}^{\delta^{_G}}\not=\! 0$ and the set of points of the $k\-$algebra $(k_*\hat Z)^\circ\,\hat b_\delta$
where  we denote by ${\rm \bar Br}_Q (b)$  the image of  ${\rm Br}_Q (b)$ in $k\bar C_G (Q)$ and by 
 $\hat b_\delta$ the image  of ${\rm \bar Br}_Q (b){\rm Tr}_{\bar C_G (Q_\delta)}^{\bar C_G (Q)}
 (\bar f^{^H})$ in the right-hand member of isomorphism~£3.6.5.
 \eject
 
 \bigskip
 \noindent
 {\bf Proposition~£3.7}\phantom{.} {\it  With the the notation above, the idempotent $\,\hat b_\delta$ is primitive 
 in~$Z(k_*\hat Z)^{E_G (Q_\delta)}\,.$ In particular, if $E_G (Q_\delta)$
 acts trivially on $\hat Z$ then $P_{\gamma^{^G}}$ contains $Q_{\delta^{^G}}$ for any local point $\delta^{^G}\!\!$ of $Q$ on~$\O Gb$ such that 
 ${\rm m}_\delta^{\delta^{^G}}\!\!\not= 0\,.$\/}
 \medskip 
 \noindent
 {\bf Proof:} Since $Q = H\cap P\,,$ for any $a\in (\O G)^P$ it is easily checked that
$${\rm Br}_Q\big({\rm Tr}_P^G(a)\big) = {\rm Tr}_P^{N_G (Q)}\big({\rm Br}_Q (a)\big)
\eqno£3.7.1\phantom{.}$$
and, in particular, we have ${\rm Br}_Q\big((\O G)_P^G\big) \cong kC_G(Q)_P^{N_G (Q)}$
(cf.~£2.9.1); consequently, since the idempotent $b\in (\O G)_P^G$ is primitive in $Z(\O G)\,,$
setting $E_G (Q) = N_G (Q)/Q\.C_G(Q)\,,$
${\rm Br}_Q(b)$ is still primitive in [17, Proposition~3.23]
 $$kC_G(Q)^{N_G (Q)} = Z\big(kC_G (Q)\big)^{E_G (Q)}
 \eqno£3.7.2,$$
 which amounts to saying that $N_G (Q)$ acts transitively over the set of $k\-$blocks of $C_G (Q)$ involved in 
 ${\rm Br}_Q(b)\,;$ hence, since any $k\-$block of $C_G (Q)$ maps on a $k\-$block of $\bar C_G (Q)$ (cf.~£2.9),
 ${\rm \bar Br}_Q (b)$ is also primitive in  $Z\big(k\bar C_G (Q)\big)^{E_G (Q)}$
  and then, it suffices to apply isomorphism~£3.6.5.

  \smallskip
On the other hand, identifying $(\O G)(Q)$ with $kC_G (Q)$ (cf.~£2.9.1), it is easily checked that ${\rm Br}_Q 
\big((\O G)^P\big) = kC_G (Q)^P$ and therefore, for any  $i\in \gamma^{_G},$ the idempotent ${\rm Br}_Q (i)$ is primitive in 
$kC_G (Q)^P$ [17, Proposition~3.23]; thus, since the canonical $P\-$algebra homomorphism 
$kC_G(Q)\to k\bar C_G (Q)$ is a {\it strict semicovering\/}
[16,~Theorem~2.9], it follows from [6, Proposition~3.15] that the image ${\rm \bar Br}_Q (i)$ of ${\rm Br}_Q (i)$ 
in $k\bar C_G (Q)^P$ remains a primitive idempotent and that, denoting by $\bar\gamma^{_G}$ the point of $P$ 
on $k\bar C_G (Q)$ determined by ${\rm \bar Br}_Q (i)\,,$ $P_{\bar\gamma^{_G}}$~remains a maximal local pointed group 
 on  the $N_G (Q)\-$algebra $k\bar C_G (Q)\,.$

\smallskip
 Moreover, since $P$ fixes~$f^{^H}\!\!$ (cf.~£3.3),  we may choose $i\in \gamma^{_G}$
fulfilling ${\rm Br}_Q (i) = {\rm Br}_Q (i) f^{^H}\!\!\,;$ in this case, it follows from isomorphism~£3.5.2 and
from [19, Proposition~3.5] that ${\rm \bar Br}_Q (i)$ is a primitive idempotent of
$\big(k\bar C_G (Q_\delta)\bar f^{^H}\big)^P$ and that $P_{\bar\gamma^{_G}}$ is also a maximal local pointed group  
on the $N_G (Q_\delta)\-$alge-bra~$k\bar C_G (Q_\delta)\bar f^{^H}\!\!.$

\smallskip
 But, it follows from isomorphism~£3.6.4 that we have
$$\big(k\bar C_G (Q_\delta)\bar f^{^H}\big)(P)\cong \big(k\bar C_H (Q)\bar f^{^H}\big)(P)\otimes_k 
(k_*\hat Z)^\circ (P)
\eqno £3.7.3\phantom{.}$$
and therefore, since evidently $ib = i\,,$ $P_{\bar\gamma^{_G}}$ determines a maximal local pointed group  
$P_{\skew1\hat{\bar\gamma}^{_G}}$ on $(k_*\hat Z)^{\circ}\,\hat b_\delta$
[9, Theorem~£5.3 and~Proposition~5.9]. Moreover,  if $E_G (Q_\delta)$
 acts trivially on $\hat Z$ then $\hat b_\delta$ is a block of $\hat Z$
and  all the  maximal local pointed groups on the $N_G (Q_\delta)\-$algebra $(k_*\hat Z)^{\circ}\,\hat b_\delta$ are mutually 
$N_G (Q_\delta)\-$con-jugate (cf.~£2.5); in this case, since $N_G (Q_\delta)$ acts trivially on the set of points of  
the $k\-$algebra $(k_*\hat Z)^{\circ}\,\hat b_\delta\,,$ $P_{\skew1\hat{\bar\gamma}^{_G}}$ contains 
$Q_{\skew2\hat{\bar\delta}^{^G}}$ for any local point $\skew2\hat{\bar\delta}^{^G}$ of $Q$ on 
$(k_*\hat Z)^{\circ}\,\hat b_\delta\,,$ which amounts to saying that any idempotent 
$\skew1\hat{\bar \imath}\in \skew1\hat{\bar\gamma}^{_G}$ has a nontrivial image in all the simple quotients of 
$(k_*\hat Z)^{\circ}\,\hat b_\delta$ (cf.~2.2.2), and   the last statement follows from~£3.6.

 \bigskip
 \noindent
 {\bf Proposition~£3.8}\phantom{.} {\it Let  $\delta^{_G}$ be a local point of $Q$ on $\O G$ such that 
${\rm m}_\delta^{\delta^{^G}}\!\!\not= 0\,.$ The commutator in $\hat N_G (Q_\delta)/Q\.C_H (Q)$
induces a group homomorphism
$$\varpi : F\too {\rm Hom} (Z, k^*)
\eqno £3.8.1\phantom{.}$$
and ${\rm Ker}(\varpi)$ is contained in  $ E_H(Q_{\delta^{^G}})\,.$ In particular, $ E_H(Q_{\delta^{^G}})$ is normal 
in~$F\,,$ $F/E_H(Q_{\delta^{^G}})$ is an Abelian $p'\-$group and, denoting by $\hat K^\delta\!$ 
and $\hat K^{\delta^{^G}}\!$ the respective converse images in $\hat C_G (Q_\delta)$ of the fixed points of
$F$ and $E_H (Q_{\delta^{^G}})$ over $\hat Z\,,$ we have the exact sequence
$$1\too \hat K^\delta\too \hat K^{\delta^{^G}}\too {\rm Hom}\big(F/E_H(Q_{\delta^{^G}}), 
k^*\big)\too 1Ê
\eqno £3.8.2.$$\/}

\par
\noindent
{\bf Proof:} It is quite clear that $F$ and $Z$ are normal subgroups of the quotient
$N_G (Q_\delta)/Q\.C_H (Q)$ and therefore their converse images $\hat F$ and $\hat Z$ in~ the quotient 
$\hat N_G (Q_\delta)/Q\.C_H (Q)$ still normalize each other; but, since we have
$$N_H(Q_\delta)\cap C_G (Q_\delta) = C_H (Q)\
\eqno £3.8.3,$$
 their commutator is contained in~$k^*\,;$ hence, indentifying ${\rm Hom} (Z, k^*)$ with the group of the automorphisms of  the $k^*\-$group $\hat Z$ which act trivially on $Z\,,$ 
we easily get homomorphism~£3.8.1.

\smallskip
In particular, ${\rm Ker}(\varpi)$ acts trivially on the $k^*\-$group $\hat Z$ 
and therefore, since its action is compatible with the bijection in~£3.6 above, it is contained in 
$ E_H(Q_{\delta^{^G}})\,;$ hence, since the $p'\-$group ${\rm Hom}(Z, k^*)$ is Abelian,
 $ E_H(Q_{\delta^{^G}})$ is normal in $ E_H(Q_\delta)$ (cf.~£3.5.3) and
$F/E_H(Q_{\delta^{^G}})$ is Abelian.

\smallskip
Symmetrically, the commutator in $\hat N_G (Q_\delta)/Q\.C_H (Q)$ also induces surjective group homomorphisms
$$\eqalign{ \hat C_G (Q_\delta)&\too {\rm Hom}\big(F/{\rm Ker}(\varpi), k^*\big)\cr
\hat C_G (Q_\delta)&\too {\rm Hom}\big(E_H(Q_{\delta^{^G}})/{\rm Ker}(\varpi), k^*\big)\cr}
\eqno £3.8.4\phantom{.}$$
and it is quite clear that the kernels respectively coincide with $\hat K^\delta$ and $\hat K^{\delta^{^G}}\!\!\,;$
consequently, the kernel of the surjective group homomorphism 
$$\hat C_G (Q_\delta)/\hat K^\delta\too \hat C_G (Q_\delta)/\hat K^{\delta^{^G}}
\eqno £3.8.5\phantom{.}$$
is canonically isomorphic to ${\rm Hom}\big(F/E_H(Q_{\delta^{^G}}),  k^*\big)\,.$ We are done.
\eject

\medskip
£3.9\phantom{.} Assume that $b$ is an inertial block of $G$ or, equivalently, that there is a primitive Dade 
$P\-$algebra $S$ such that, with the  notation in~£2.13 above, we have a $P\-$interior algebra isomorphism
$$(\O G)_{\gamma^{_G}}\cong S\otimes_\O \O_*\hat L
\eqno £3.9.1\phantom{.}$$
where we consider $S$ endowed with the unique $P\-$interior algebra structure fulfilling ${\rm det}_S (P) = \{1\}$ 
(cf.~£2.14). In this case,  it follows from [6, Lemma~1.17] and [8, proposition~2.14 and Theorem~3.1] that (cf.~£2.8)
$$E = F_{\O G}(P_{\gamma^{_G}}) =
F_S(P_{\{1_S\}}) \cap F_{\O_*\hat L}(P_{\{1_{\hat L}\}})
\eqno £3.9.2\phantom{.}$$
and, in particular, that $S$ is {\it $E\-$stable\/} [8, Proposition~2.18]. Moreover, since we have a $P\-$interior
algebra embedding (cf.~£2.14)
$$\O\too {\rm End}_\O (S) \cong S^\circ\otimes_\O S
\eqno £3.9.3,$$
we still have a $P\-$interior algebra embedding 
$$\O_*\hat L\too S^\circ\otimes_\O (\O G)_{\gamma^{_G}}
\eqno £3.9.4.$$

\medskip
£3.10\phantom{.} Conversely, always  with the notation in~£2.13, assume that $S$ is an $E\-$stable  Dade $P\-$algebra
or, equivalently, that $E$ is contained in $F_S (P_\pi)$ where $\pi$ denotes the unique local point of $P$ on $S$ 
(cf.~£2.14); since we have [9,~Proposition~5.9]
$$F_S (P_\pi)\cap F_{\O G}(P_{\gamma^{_G}})\i F_{S^\circ\otimes_\O \O G}
(P_{\pi\times \gamma^{_G}})
\eqno £3.10.1$$
where $\pi\times \gamma^{_G}$ denotes the local point of $P$ on $S^\circ\otimes_\O \O G$
determined by $\pi$ and $\gamma^{_G}$ [9, Proposition~5.6], and we still have [18, Theorem~9.21]
$$\hat F_S (P_\pi)\cong k^*\times F_S (P_\pi)
\eqno £3.10.2,$$
it follows from [9, proposition~5.11] that the $k^*\-$group $\hat E$ is isomorphic to a $k^*\-$subgroup of 
$\hat F_{S^\circ\otimes_\O \O G}(P_{\pi\times \gamma^{_G}})\,;$ then, since $E$ 
is a $p'\-$group, it follows from [10, Proposition~7.4] that there is an injective unitary $P\-$interior algebra homomorphism
$$\O_*\hat L\too  (S^\circ\otimes_\O \O G)_{\pi\times \gamma^{_G}}
\eqno £3.10.3\phantom{.}$$ 
and, in particular, we have 
$$\vert P\vert\vert E\vert\le {\rm rank}_\O (S^\circ\otimes_\O \O G)_{\pi\times \gamma^{_G}}
\eqno £3.10.4.$$

\bigskip
\noindent
{\bf Proposition~£3.11}\phantom{.} {\it With the notation above, the block $b$ is inertial if and only
if there is  an $E\-$stable  Dade $P\-$algebra $S$ such that 
$${\rm rank}_\O \big(S^\circ\otimes_\O \O G\big)_{\pi\times \gamma^{_G}}
= \vert P\vert\vert E\vert
\eqno £3.11.1$$\/}

\par
\noindent
{\bf Proof:} If $b$ is inertial then the equality~£3.11.1 follows from the existence of embedding~£3.9.4.
\eject

\smallskip
Conversely, we claim that if equality~£3.11.1 holds then the corresponding homomorphism~£3.10.3
is an isomorphism; indeed, since this homomorphism is injective and we have ${\rm rank}_\O (\O_*\hat L) = 
\vert P\vert\vert E\vert\,,$ it suffices to prove that the reduction to $k$ of homomorphism~£3.10.3
remains injective;  but, according to [10, 2.1], it also follows from [10, Proposition~7.4] that, setting 
${}^k S = k\otimes_\O S\,,$ there is an injective unitary $P\-$interior algebra homomorphism
$$k_*\hat L\too ({}^k S^\circ\otimes_k k G)_{\bar\pi\times \bar\gamma^{_G}}
\eqno £3.11.2,$$ 
where $\bar\pi$ and $\bar\gamma^{_G}$ denote the respective images of $\pi $ and $\gamma^{_G}$ in
${}^k S^\circ$ and $k G\,,$ which is a conjugate of the reduction to $k$ of homomorphism~£3.10.3.

\smallskip
Now, embedding~£3.9.3 and the structural embedding
$$(S^\circ\otimes_\O \O G)_{\pi\times \gamma^{_G}}\too S^\circ\otimes_\O 
(\O G)_{\gamma^{_G}}
\eqno £3.11.3\phantom{.}$$
determine  $P\-$interior algebra embeddings
$$\matrix{ S\otimes_\O (S^\circ\otimes_\O \O G)_{\pi\times \gamma^{_G}}&\too 
&S\otimes_\O S^\circ\otimes_\O(\O G)_{\gamma^{_G}}\cr
\hskip-20pt \wr\Vert&\phantom{\Big\uparrow}&\big\uparrow\cr
\hskip-20pt S\otimes_\O \O_*\hat L&&(\O G)_{\gamma^{_G}}\cr}
\eqno £3.11.4;$$
thus, since $P$ has a unique local point on $S\otimes S^\circ\otimes_\O(\O G)_{\gamma^{_G}}$ [9, Theorem~5.3],
we get a $P\-$interior algebra embedding
$$(\O G)_{\gamma^{_G}}\too S\otimes_\O \O_*\hat L
\eqno £3.11.5$$
which proves that $b$ is inertial. We are done.

\medskip
£3.12\phantom{.} With the notation above, assume that the block $b$ is inertial; then, denoting by $\chi$ the unique
local point of $Q$ on $S$ (cf.~£2.14) and by $\delta^{^G}$ a local point of $Q$ on $\O Gb$ such that
${\rm m}_\delta^{\delta^{^G}}\not=0\,,$ there is a unique local point $\hat\delta^{^L}$ of $Q$ on $\O_*\hat L$ 
such that isomorphism~£3.9.1 induces  a $Q\-$interior algebra embedding  [9, Proposition~5.6]
$$(\O G)_{\delta^{^G}}\too S_\chi\otimes_\O (\O_*\hat L)_{\hat\delta^{^L}}
\eqno £3.12.1;$$
but, the image of $Q$ in $(S_\chi)^*$ need not be contained in the kernel of the corresponding {\it determinant map\/}.
Note that,  as above, it follows from this embedding and  from [6,~Lemma~1.17] and [8, Proposition~2.14 and Theorem~3.1] that
$$E_G (Q_{\delta^{^G}}) = F_{\O G}(Q_{\delta^{^G}}) = F_{S}(Q_\chi)\cap
F_{\O_*\hat L}(Q_{\hat\delta^{^L}})
\eqno £3.12.2,$$
so that {\it the Dade $Q\-$algebra $S_\chi$ is $E_G (Q_{\delta^{^G}})\-$stable\/}; as in~£2.13 above, let us
consider the corresponding semidirect products
$$M = Q\rtimes F\qq \hat M = Q\rtimes \hat F
\eqno £3.12.3.$$
We are ready to state our main result.
\eject

\bigskip
\noindent
{\bf Theorem~£3.13}\phantom{.} {\it With the notation above, assume that the block $b$ of $G$ is inertial.
Then, there is a $Q\-$interior algebra isomorphism
$$(\O H)_\delta\cong S_\chi\otimes_\O \O_*\hat M 
\eqno £3.13.1\phantom{.}$$
and, in particular, the block $c$ of $H$ is inertial too.\/}

\medskip
\noindent
{\bf Proof:} We argue by induction on $\vert G/H\vert\,;$ in particular, if $H'$ is a proper normal subgroup of $G$
which properly contains $H\,,$ it suffices to choose a block $c'$ of~$H'$ fulfilling $c'b\not= 0$ to get $c'c\not= 0$
and the induction hypothesis successively proves that the block $c'$ of $H'$  is inertial and then that the block $c$
is inertial too; moreover, setting $Q' = H'\cap P\,,$ the corresponding Dade $Q'\-$algebra comes from $S$
and therefore the final Dade $Q\-$algebra also comes from $S\,.$ Consequently, since $G$ fixes $c\,,$ it follows from the {\it Frattini argument\/} that we have (cf.~£2.3)
 $$G = H\. N_G (Q_\delta)
 \eqno £3.13.2$$
  and therefore we may assume that either  $C_G (Q_\delta)\i H$ or  
$G = H\. C_G (Q_\delta)\,.$

\smallskip
Firstly assume that $C_G (Q_\delta)\i H\,;$ in this case, it follows from [18, Proposition~15.10] that $b = c\,;$
moreover, since $C_G (Q_\delta) = C_H (Q)\,,$ it follows from~£3.6 above that $Q$ has a unique local point
$\delta^{^G}\!\!$ on $\O Gb$ such that ${\rm m}_\delta^{\delta^{^G}}\!\! \not= 0\,,$ and from isomorphism~£3.6.4 that we have
$$(\O H)(Q_\delta)\cong k\bar C_H (Q)\bar f^{^H}\cong k\bar C_G (Q_\delta)\bar f^{^H}
\eqno £3.13.3;$$
in particular, $N_G (Q_\delta)$ normalizes $Q_{\delta^{^G}}$
and therefore the inclusion~£3.5.3 becomes an equality
$$E_G (Q_{\delta^{^G}}) = E_G (Q_\delta) 
\eqno £3.13.4;$$ 
thus, since $F$
is obviously contained in $E_G (Q_\delta)\,,$ $S_\chi$ is $F\-$stable too.
Consequently, according to Proposition~£3.11, it suffices to prove that
$${\rm rank}_\O (S^\circ_\chi\otimes_\O \O H)_{\chi\times \delta}
= \vert Q\vert\vert F\vert
\eqno £3.13.5.$$

\smallskip
As in~£3.12 above, the $P\-$interior algebra embedding~£3.9.4 induces a $Q\-$interior algebra embedding
[9, Theorem~5.3]
$$(\O_*\hat L)_{\hat\delta^{^L}}\too S_\chi^\circ\otimes_\O (\O G)_{\delta^{_G}}
\eqno £3.13.6\phantom{.}$$
and it suffices to apply again  [6, Lemma~1.17] and [8, Proposition~2.14 and Theorem~3.1] to get
$$E_L (Q_{\hat\delta^{^L}}) = F_{\O_*\hat L}(Q_{\hat\delta^{^L}}) = 
F_{S}(Q_\chi)\cap F_{\O G}(Q_{\delta^{^G}})
\eqno £3.13.7,$$
so that we obtain
$$E_L (Q_{\hat\delta^{^L}}) = E_G (Q_{\delta^{^G}})\i F_{S}(Q_\chi)
\eqno £3.13.8.$$
In particular, it follows from  [8, Proposition~2.12] that for any $x\in N_G (Q_\delta)$ there is $s_x\in (S_\chi)^*$ fulfilling 
$$s_x\.u = u^x\. s_x
\eqno £3.13.9\phantom{.}$$
for any $u\in Q\,,$ and therefore, choosing a set of representatives $X\i N_G(Q_\delta)$ for $G/H$ (cf.~£3.13.2),
we get an $\O Q\-$bimodule direct sum decomposition
$$S_\chi^\circ \otimes_\O \O G = \bigoplus_{x\in X} (s_x\otimes x)(S_\chi^\circ \otimes_\O \O H)
\eqno £3.13.10.$$

\smallskip
But, for any $x\in N_G (Q_\delta)\,,$ the element $s_x\otimes x$ normalizes the image of~$Q$ in 
$S_\chi^\circ \otimes_\O \O H$ and it is clear that it also normalizes the local point $\chi\times \delta$
of~$Q$ on this $Q\-$interior algebra; more precisely, if $S_\chi = \ell S\ell$ for $\ell\in \chi$
and $(\O H)_\delta = j(\O H) j$ for $j\in \delta\,,$ there is $j'\in \chi\times\delta$ such that [9, Proposition~5.6]
$$j'(\ell\otimes j) = j' = (\ell\otimes j)j'
\eqno £3.13.11;$$
thus, for any $x\in N_G (Q_\delta)$ the idempotent $j'^{s_x\otimes x}$ still belongs to $\chi\times \delta$
and therefore there is an inversible element $a_x$ in $(S_\chi^\circ \otimes_\O \O H)^Q$ fulfilling
$$j'^{s_x\otimes x} = j'^{a_x}
\eqno £3.13.12,$$
so that we get the new $\O Q\-$bimodule direct sum decomposition
$$j'(S^\circ \otimes_\O \O G)j' = \bigoplus_{x\in X} (s_x\otimes x)(a_x)^{-1}j'(S^\circ \otimes_\O \O H)j'
\eqno £3.13.13.$$

\smallskip
Moreover, the equality in~£3.13.8 forces the group $E_G (Q_\delta) = E_G (Q_{\delta^{^G}})$ to have a normal 
Sylow $p\-$subgroup  and therefore, since we are assuming that $C_G (Q_\delta)\i H\,,$ it follows from equality~£3.13.2 that the quotient $G/H$ also has a normal Sylow $p\-$subgroup. At this point, arguing by induction, we may assume that $G/H$ is either a $p\-$group or a $p'-$group.

\smallskip
Firstly assume that $G/H$ is a $p\-$group or, equivalently, that $G = H\. P$ [9, Lemma~3.10]; in this case,
it follows from [6, Proposition~6.2] that the inclusion homomorphism $\O H\to \O G$ is a {\it strict
semicovering\/} of $Q\-$interior algebras (cf.~£2.1) and, in particular, we have $\delta\i \delta^{^G}$ since  
${\rm m}_\delta^{\delta^{^G}}\not= 0\,;$ similarly, since for any subgroup $R$ of $Q$ we have [9, Proposition~5.6]
$$\eqalign{(S^\circ\otimes_\O \O H)(R) &\cong S(R)^\circ\otimes_k (\O H)(R)\cr
(S^\circ\otimes_\O \O G)(R) &\cong S(R)^\circ\otimes_k (\O G)(R)\cr}
\eqno £3.13.14,$$
\eject
\noindent
it follows from  [6, Theorem~3.16] that the corresponding $Q\-$interior algebra homomorphism
$S^\circ\otimes_\O\O H\to S^\circ\otimes_\O\O G$ is also a {\it strict semicovering\/} and, in particular, we have 
$\chi\times \delta\i \chi\times\delta^{^G}\,,$ so that $j'$ belongs to $\chi\times\delta^{^G}\,.$

\smallskip
But, since $Q_{\delta^{^G}}\i P_{\gamma^{^G}}$ (cf.~£3.4), it is easily checked that $Q_{\chi\times\delta^{^G}}\i 
P_{\pi\times\gamma^{^G}}\,,$ where as above $\pi$ is the unique local point of $P$ on $S\,,$ and therefore we get
 the  $Q\-$interior algebra embedding (cf.~embeddings~£2.2.3 and~£3.9.4)
 $$(S^\circ\otimes_\O \O G)_{\chi\times \delta^{^G}}\too {\rm Res}_Q^P(S^\circ\otimes_\O \O G)_{\pi\times \gamma^{^G}}  \cong {\rm Res}_Q^P(\O_*\hat L)
 \eqno £3.13.15;$$
in particular, it follows from equality~£3.13.13 that we have
$$\vert X\vert\, {\rm rank}_\O\, (S^\circ_\chi\otimes_\O \O H)_{\chi\times \delta}\le \vert L\vert
\eqno £3.13.16.$$
Moreover, we have $\vert X\vert = \vert G/H\vert = \vert P/Q\vert$ and, since $C_P (Q)\i Q\,,$ it follows from
[4, Ch. 5, Theorem~3.4] that a subgroup of $E\i L$ which centralizes $Q = H\cap P$ still centralizes $P\,,$
so that $E$ acts faithfully on $Q\,;$ in particular, $\hat\delta^{^L}$ is the unique local point of $Q$ on $\O_*\hat L$ (actually, we have $\hat\delta^{^L} = \{1_{_{\O_*\hat L}}\})$ and therefore, since 
(cf.~£3.13.4 and~£3.13.8)
$$E_L (Q_{\hat\delta^{^L}}) = E_G (Q_{\delta^{^G}}) = E_G (Q_\delta)\j F
\eqno £3.13.17\phantom{.}$$
and $E_G (Q_\delta)/F$ is a $p\-$group, the $p'\-$group $E$ is actually isomorphic to $F\,.$

\smallskip
 Consequently, it follows from the inequalities £3.10.4 and £3.13.16 that
$$\vert F\vert \vert Q\vert \le {\rm rank}_\O\, (S^\circ_\chi\otimes_\O \O H)_{\chi\times \delta}
\le  \vert L\vert/\vert X\vert = \vert F\vert \vert Q\vert 
\eqno £3.13.18\phantom{.}$$
which forces equality~3.13.6.

\smallskip
Secondly assume that $G/H$ is a $p'\-$group; in this case, we have $Q = P\,,$ $\delta = \gamma$ and
$\delta^{^G} = \gamma^{^G}\,;$ in particular, since we are assuming that 
$$C_G (Q_\delta)\i H\qq E_G (Q_{\delta^{^G}})  = E_G (Q_\delta)
\eqno £3.13.19,$$
we actually get 
$$\vert X\vert = \vert G/H\vert = \vert E_G (P_{\gamma^{^G}})\vert/\vert E_H (Q_\delta)\vert 
= \vert E\vert/\vert F\vert
\eqno £3.13.20.$$
Moreover, we claim that, as above, the idempotent $j'$ remains primitive in~$(S\otimes_\O \O G)^P$
{\footnote{\dag}{\cds The corresponding argument has been forgotten in [18] at the end of the proof of 
Proposition~15.19!}}, so that 
it belongs to $\pi\times \gamma^{_G}\,;$ indeed, setting
$$A' = j'(S^\circ \otimes_\O \O G)j'\qq B '= j'(S^\circ \otimes_\O \O H)j'
\eqno £3.13.21,$$ let $i'$ be a primitive idempotent of $A'^P$ such that ${\rm Br}_P(i')\not= 0\,;$
in particular, $i'$ belongs to $\pi\times \gamma^{_G}$ and we may assume that
$$i'A'i' = (S^\circ\otimes_\O \O G)_{\pi\times \gamma^{^G}}  \cong \O_*\hat L
\eqno £3.13.22.$$
\eject

\smallskip
It is clear that the multiplication by $B'$ on the left and the action of~$P$ by conjugation endows $A'$ with
a $B'P\-$module structure and, since the idempotent $j'$ is primitive in $B'^P\,,$ equality~£3.13.13 provides 
a direct sum decomposition of $A'$ in indecomposable $B'P\-$modules. More explicitly, note that $B'$ is an indecomposable $B'P\-$module since we have ${\rm End}_{B'P}(B') = B'^P\,;$ but,
for any $x\in X\,,$ the inversible element 
$$a'_x = (s_x\otimes x)(a_x)^{-1}j'
\eqno £3.13.23\phantom{.}$$
 of $A'$ together with the action of $x$ on $P$ determine an automorphism $g_x$ of~$B'P\,;$ thus,  
 equality~£3.13.13 provides the following direct sum decomposition on indecomposable $B'P\-$modules
 $$A'\cong \bigoplus_{x\in  X} {\rm Res}_{g_x}(B')
 \eqno £3.13.24.$$

 \smallskip
 Moreover, we claim that the $B'P\-$modules ${\rm Res}_{g_x}(B')$ and ${\rm Res}_{g_{x'}}(B')$ for $x,x'\in X$
 are isomorphic if and only if $x = x'\,;$ indeed, a $B'P\-$module isomorphism
 $${\rm Res}_{g_x}(B')\cong {\rm Res}_{g_{x'}}(B')
 \eqno £3.13.25\phantom{.}$$
 is necessarily determined by the multiplication on the right by an inversible element $b'$ of $B'$ fulfilling
 $$(xux^{-1})\.b' = b'\.(x'ux'^{-1})
 \eqno £3.13.26\phantom{.}$$
or, equivalently, $(u\.j')^{b'} = u^{xx'^{-1}}\.j'$ for any $u\in P\,,$ which amounts to saying that the automorphism of $P$ determined by the conjugation by $x'x^{-1}$ is a {\it $B'\-$fusion\/} from $P_\gamma$ to $P_\gamma$ [8,~Proposition~2.12]; but, once again from  [6,~Lemma~1.17] and [8, Proposition~2.14 and Theorem~3.1] we get
$$F_{A'} (P_{\gamma^{^G}})= E_G (P_{\gamma^{^G}}) = E\qq F_{B'}(P_\gamma) = E_H (P_\gamma)
\eqno £3.13.27\,;$$
hence our claim now follows from qualities~£3.13.20.

\smallskip
On the other hand, it is clear that $A'i'$ is a direct summand of the $B'P\-$module $A'$ and therefore there is 
$x\in X$ such that ${\rm Res}_{g_x}(B')$ is a direct summand of the $B'P\-$module $A'i'\,;$ but, it follows from
[8, Proposition~2.14] that we have
$$F_{i'A'i'} (P_{\gamma^{^G}}) = F_{A'} (P_{\gamma^{^G}}) = E
\eqno £3.13.28\phantom{.}$$
and therefore, once again applying [8,~Proposition~2.12], for any $y\in N_G(P_{\gamma^{^G}})$
there is an inversible element $c'_y$ in $A'$ fulfilling  
$$c'_y(u\.i')(c'_y)^{-1} = yuy^{-1}\.i'
\eqno £3.13.29\phantom{.}$$
for any $u\in P\,;$ then, for any $x'\in X\,,$ it is clear that $A'i' = A'i'c'_{x^{-1}x'}$ has a direct summand
isomorphic to  ${\rm Res}_{g_{x'}}(B')\,,$ which forces the equality of the $\O\-$ranks of $A'i'$ and $A'\,,$
so that $A'i' = A'$ and $i' =j'\,,$ which proves our claim. Consequently, it follows from the equalities £3.13.13 and £3.13.20 that
$${\rm rank}_\O\, (S^\circ_\chi\otimes_\O \O H)_{\chi\times \delta}
=  \vert L\vert/\vert X\vert = \vert F\vert \vert Q\vert 
\eqno £3.13.30,$$
so that equality~3.13.6 holds.

\smallskip
From now on, we assume that $H\.C_G(Q_\delta) = G\,;$ in particular, $C_G (Q)$ stabilizes $\delta\,,$ we have
$E_G(Q_\delta) = E_H (Q_\delta) = F$ and we can choose the set of representatives $X$ for $G/H$ contained in 
$C_G (Q)\,,$ so that this time we get the $\O Q\-$bimodule direct sum decomposition
$$S_\chi^\circ \otimes_\O \O G = \bigoplus_{x\in X} (1_{_S}\otimes x)(S_\chi^\circ \otimes_\O \O H)
\eqno £3.13.31.$$
Since any $z\in C_G (Q)$ stabilizes $\delta\,,$  choosing again $\ell\in \chi\,,$ $j\in \delta$ and $j'\in \chi\times\delta$ such that [9, Proposition~5.6]
$$j'(\ell\otimes j) = j' = (\ell\otimes j)j'
\eqno £3.13.32,$$
there is  an inversible element $a_z$ in $(\O H)^Q$ fulfilling $j^{z} = j^{a_z}\,;$
consequently, with the notation above, from these choices and equality~£3.13.31 we have
$$A' = \bigoplus_{x\in X} (1_{_S}\otimes x(a_x)^{-1})B'
\eqno £3.13.33.$$
 As in Proposition~£3.8, denote by $\hat K^\delta$ the converse image in $\hat C_G (Q)$ of the fixed points of $F$
 in $\hat Z$ and by $K^\delta$ the {\it $k^*\-$quotient\/} 
$\hat K^\delta/k^*$ of~$\hat K^\delta\,;$ since $\hat K^\delta$ is a normal 
$k^*\-$subgroup of $\hat C_G (Q)\,,$ $H\.K^\delta$ is a normal subgroup of~$G$
and therefore, arguing by induction, we may assume that it coincides with~$H$ or with~$G\,.$

\smallskip
Firstly assume that $H\.K^\delta = G\,;$ in this case, since we have $K^\delta = C_G(Q)\,,$  $F$ acts trivially on 
$\hat Z$ and we have  $F = E_H (Q_{\delta^{^G}})$ for any local point $\delta^{^G}$ of $Q$ on $\O Gb$ such that
${\rm m}_\delta^{\delta^{^G}}\not= 0\,,$ so that $S_\chi$ is $F\-$stable (cf.~£3.12.2); consequently, 
according to Proposition~£3.11, once again it suffices to prove that
$${\rm rank}_\O (S^\circ_\chi\otimes_\O \O H)_{\chi\times \delta} = \vert Q\vert\vert F\vert
\eqno £3.13.34.$$

\smallskip
For any $z\in C_G(Q)\,,$ the element $z(a_z)^{-1}$
stabilizes $j(\O H)j = (\O H)_\delta$ and actually it induces a $Q\-$interior algebra automorphism $g_z$ of the source algebra $(\O H)_\delta\,;$ but, symmetrically, $C_G (Q)$ acts trivially on [8, Proposition~2.14 and~Theorem~3.1]
$$\hat F = \hat E_H (Q_\delta)^\circ \cong \hat F_{(\O H)_\delta} (Q_\delta)\eqno £3.13.35;$$
hence, it follows from [10, Proposition~14.9] that $g_z$ is an {\it inner automorphism\/} and therefore, up to modifying
our choice of $a_z\,,$ we may assume that $z(a_z)^{-1}$ centralizes $(\O H)_\delta\,;$ then, for any $x\in X$ the element $1_{_S}\otimes x(a_x)^{-1}$ centralizes 
$$ B '= j'(S^\circ \otimes_\O \O H)j'
\eqno £3.13.36.$$
and therefore, denoting by $C$ the centralizer of $B'$ in $A'\,,$ it follows from equality~£3.13.33 that we have
$$A' = C\otimes_{Z(B')} B'
\eqno £3.13.37;$$
in particular, we get $A'^Q = C\otimes_{Z(B')} B'^Q$ which induces a $k\-$algebra isomorphism [10, 14.5.1]
$$A' (Q)\cong C\otimes_{Z(B')} kZ(Q)
\eqno £3.13.38\phantom{.}$$
and then it follows from isomorphism~£3.6.4 that
$$k\otimes_{Z(B')} C\cong (k_*\hat Z)^\circ
\eqno £3.13.39.$$

\smallskip
At this point,  for any local point $\delta^{^G}$ of $Q$ on~$\O Gb$ such that ${\rm m}_\delta^{\delta^{^G}}
\!\!\not= 0\,,$ it follows from Proposition~£3.7 that $Q_{\delta^{^G}}\i P_{\gamma^{^G}}\,,$ so that $Q_{\chi\times \delta^{^G}}\i P_{\pi\times \gamma^{^G}}$ [9,~Proposition~5.6] and therefore $\chi\times \delta^{^G}$ is also a local point 
of $Q$ on the $P\-$interior algebra (cf.~embedding~£3.9.4)
$$(S^\circ \otimes_\O \O G)_{\pi\times  \gamma^{^G}}\cong \O_*\hat L
\eqno £3.13.40;$$
actually, since $N_G (P)$ normalizes $Q = H\cap P\,,$ $Q$ is normal in $L$ and therefore all the points of $Q$ on
$\O_*\hat L$ are local (cf.~£2.10). In conclusion, since $\{1_{_L}\}$ is the  unique point of  $P$ on $\O_*\hat L\,,$
isomorphism~£3.13.40 induces a bijective correspondence between the sets of local points of $Q$ on 
$$j'(S^\circ \otimes_\O \O Gb)j' = A'(1\otimes b)
\eqno £3.13.41\phantom{.}$$
 and on~$\O_*\hat L\,;$ moreover, note that if two local points
 $\chi\times\delta^{^G}$ and $\chi\times\varepsilon^{_G}$ of~$Q$ on the left-hand member of~£3.13.40 correspond to  two local points  $\hat\delta^{^G}$ and $\hat\varepsilon^{_G}$ of $Q$ on~$\O_*\hat L\,,$ choosing suitable
$j^{^G}\in \delta^{^G}\,,$ $k^{^G}\in \varepsilon^{_G}\,,$  $\hat\jmath^{^G}\in \hat\delta^{^G}\!$ and
$\hat k^{^G}\in \hat \varepsilon^{_G}\,,$ from isomorphism~£3.13.40 we still get an $\O Q\-$bimodule
isomorphism
    $$j^{^G}\!\!A'k^{^G}\cong
\hat\jmath^{^G}\!(\O_*\hat L)\hat k^{^G}
  \eqno £3.13.42.$$

 \smallskip
Consequently, since we have $A'^Q = C\otimes_{Z(B')} B'^Q$ and $C$ is a free $Z(B')\-$ module, for suitable primitive idempotents
$\bar\jmath^{^G}$ and $\bar k^{^G}$ of $C$ we have (cf.~£3.13.37 and~£3.13.38)
$$\eqalign{{\rm dim}\big(k\otimes_{Z(B')} (\bar\jmath^{^G}\! C\bar k^{^G})\big)\,{\rm rank}_\O (B') &= 
{\rm rank}_\O \big(\hat\jmath^{^G}\! (\O_*\hat L)\hat k^{^G}\big)\cr
{\rm dim}\,\big(k\otimes_{Z(B')} (\bar\jmath^{^G}\! C\bar k^{^G})\big) &= 
{\rm rank}_{kZ(Q)}\, \big(\hat\jmath^{^G}\! (\O_*\hat L)\hat k^{^G}\big)(Q)\cr}
\eqno £3.13.43\,;$$
thus, since  the respective {\it multiplicities\/} (cf.~£2.2) of points $\hat\delta^{^G}$ and ${\rm Br}_Q^{\O_*\hat L}(\hat\delta^{^G})$ of $Q$ on $\O_*\hat L$ and on $(\O_*\hat L)(Q)\cong k_*C_{\hat L} (Q)$ coincide each other,
we finally get
$$\vert L\vert = {\rm rank}_\O (\O_*\hat L) = \vert \bar C_L (Q)\vert\,{\rm rank}_\O (B')
\eqno £3.13.44.$$

\smallskip
But, according to~£3.5.4, $N_G(P_{\gamma^{^G}})$ normalizes $\gamma$ which determines $f^{^H}$ (cf.~£3.3.1)
and therefore $\gamma$ determines the unique local point $\delta$ of $Q$ on~$\O H$ associated with $f^{^H}\,;$ thus, $N_G(P_{\gamma^{^G}})$ is contained in $N_G (Q_\delta)$ which acts tri-vially on $\hat Z\,,$
and therefore $N_G(P_{\gamma^{^G}})$ stabilizes all the local points $\delta^{^G}\!\!$ of~$Q$ on $\O Gb$ fulfilling 
${\rm m}_\delta^{\delta^{^G}}\!\!\not= 0$ (cf.~£3.6); hence, 
it follows from isomorphism~£3.13.40 above that, denoting by $\hat\delta^{^G}\!$ the  point of $Q$ 
on $\O_*\hat L$ determined by $\delta^{^G}\!\!\,,$ $L$ normalizes $Q_{\hat\delta^{^G}}\,;$ in particular, we have
$$\eqalign{F = E_G (Q_\delta) &= E_G (Q_{\delta^{^G}}) = F_{(\O G)_{\gamma^{^G}}}(Q_{\delta^{^G}})\cr
& = E_L (Q_{\hat\delta^{^G}}) = L/Q\.C_L (Q)\cr}
\eqno £3.13.45\phantom{.}$$
and therefore from equality~£3.13.44 we get
$$\vert F\vert \vert Q\vert = \vert L\vert/\vert \bar C_L (Q)\vert  = {\rm rank}_\O (B')
\eqno £3.13.46,$$
which proves that $c$ is inertial.

\smallskip
Finally, assume that $K^\delta = C_H (Q)\,;$ in this case, since the commutator in $\hat N_G (Q_\delta)/Q\.C_H (Q)$  induces a group isomorphism
$$ \hat C_G (Q_\delta)/\hat K^\delta \cong {\rm Hom}\big(F/{\rm Ker}(\varpi), k^*\big)
\eqno £3.13.47,$$ 
the quotient $G/H$ is an Abelian $p'\-$group and, in particular, we have $P = Q\,.$ But, since with our choices above
we still have (cf.~£3.13.33)
$$(\O G)_\delta = j(\O G)j = \bigoplus_{x\in X} x(a_x)^{-1}(\O H)_\delta
\eqno £3.13.48\phantom{.}$$
where the element $x(a_x)^{-1}$ determines a $Q\-$interior algebra automorphism of~$(\O H)_\delta\,,$
it suffices to consider the $k^*\-$group
$$\hat U = \bigcup_{x\in X} x(a_x)^{-1}\big((\O H)_\delta^Q\big)^*
\eqno £3.13.49\phantom{.}$$
to get the $Q\-$interior algebra $(\O G)_\delta$ as the {\it crossed product\/} [3, 1.6]
$$(\O G)_\delta \cong (\O H)_\delta \otimes_{((\O H)_\delta^Q)^*} \hat U
\eqno £3.13.50.$$
\eject

\smallskip
Then, since $G/H$ is a $p'\-$group, denoting by $U$ the $k^*\-$quotient of $\hat U$ it follows from [10, Proposition~£4.6] that the exact
sequence
$$1\too j + J\big((\O H)_\delta^Q\big)\too U\too G/H\too 1
\eqno £3.13.51\phantom{.}$$
is {\it split\/} and therefore, for a suitable central $k^*\-$extension $\widehat{ G/H}$ of $G/H\,,$ we still get an evident
 $Q\-$interior algebra isomorphism
$$(\O G)_\delta \cong (\O H)_\delta \otimes_{k^*} \widehat{ G/H}
\eqno £3.13.52;$$
at this point, it suffices to compute the {\it Brauer quotients\/} at $Q$ of both members to get
$$k\otimes_{kZ(Q)} (\O G)_\delta (Q)\cong k_*\widehat{ G/H}
\eqno £3.13.53\phantom{.}$$
and therefore, comparing this $k\-$algebra isomorphism with isomorphism~£3.6.4, we obtain a $Q\-$interior algebra isomorphism
$$(\O G)_\delta \cong (\O H)_\delta \otimes_{k^*} \hat Z^\circ
\eqno £3.13.54\phantom{.}$$
for a suitable action of $Z$ over $(\O H)_\delta$ defined, up to {\it inner automorphisms\/}
of the  $Q\-$interior algebra $ (\O H)_\delta\,,$ by the group homomorphism
$$Z\too {\rm Aut}_{k^*}\big(\hat E_{H}(Q_\delta)\big)
\eqno £3.13.55\phantom{.}$$
induced by the commutator in $\hat N_G (Q_\delta)/Q\.C_H (Q)$ [10, Proposition~14.9].

\smallskip
Similarly, considering the trivial action of $Z$ over $S\,,$ we also obtain the  $Q\-$interior algebra isomorphism
$$S^\circ \otimes_\O (\O G)_\delta \cong \big(S^\circ \otimes_\O (\O H)_\delta\big)\otimes_{k^*} \hat Z^\circ
\eqno £3.13.56;$$
since $\chi\times \delta$ is the unique local point of $Q$ on $S^\circ \otimes_\O (\O H)_\delta \,,$ 
we have $j'^{\bar z}  = j'^{b_{\bar z}}$ for a suitable inversible element $b_{\bar z}$ in $\big(S^\circ 
\otimes_\O (\O H)_\delta\big)^Q\,;$ hence, arguing as above, we finally obtain a  $Q\-$interior algebra isomorphism
$$(S^\circ \otimes_\O \O G)_{\chi\times\delta} \cong (S^\circ \otimes_\O \O H)_{\chi\times\delta}\otimes_{k^*} \hat Z^\circ
\eqno £3.13.57.$$

\smallskip
Moreover, since the $k\-$algebra $k_*\hat Z$ is now semisimple,  for any pair of primitive idempotents  $\hat\imath$ and  
$\hat\imath'$ of $\O_*\hat Z$ we have  $\hat\imath (\O_*\hat Z)\hat\imath' = \O$ {\it or\/} $\{0\}\,;$ now, since 
$\O_*\hat Z$ is a unitary $\O\-$subalgebra of $(S^\circ \otimes_\O \O G)_{\chi\times\delta}\i S^\circ \otimes_\O \O G\,,$ 
in the first case from isomorphism~£3.13.56 we get
$${\rm rank}_\O\big( \hat\imath(S^\circ \otimes_\O \O G)\hat\imath'\big)\le 
{\rm rank}_\O (S^\circ  \otimes_\O \O H)_{\chi\times\delta}
\eqno £3.13.58,$$
whereas in the second case we simply get  $\hat\imath(S^\circ \otimes_\O \O G)\hat\imath' =\{0\}\,;$
hence, since isomorphism~£3.13.57 implies that
$${\rm rank}_\O(S^\circ \otimes_\O \O G)_{\chi\times\delta} = {\rm rank}_\O(S^\circ 
\otimes_\O \O H)_{\chi\times\delta}\,\vert Z\vert\eqno £3.13.59,$$
\eject
\noindent
all the inequalities~£3.13.58 are actually equalities and, in particular, we get (cf.~embedding~£3.9.4)
$$\vert L\vert = {\rm rank}_\O(S^\circ \otimes_\O \O G)_{\pi\times\gamma^{^G}} = 
{\rm rank}_\O (S^\circ  \otimes_\O \O H)_{\chi\times\delta}
\eqno £3.13.60\phantom{.}$$
since $P = Q$ and $ \pi\times\gamma^{^G}= \chi\times \delta^{^G}$ (cf.~£3.4). Consequently, according 
to Proposition~£3.11, it suffices to prove that $S$ is $F\-$stable.

\smallskip
On the other hand, it follows from Proposition~£3.7 that $F$ acts transitively over the set of primitive
idempotents of   $Z(k_*\hat Z)\,\hat b_\delta\,;$ but, since $k_*\hat Z$ is semisimple, this set is canonically isomorphic to the set of points of this $k\-$algebra (cf.~£2.2), so that $F$ acts  transitively over the set of  local points $\delta^{^G}$ of~$Q$ on $\O Gb$ fulfilling  ${\rm m}_\delta^{\delta^{^G}}\not= 0$ (cf.~£3.6). More precisely, choosing $\delta^{^G} = \gamma^{^G}\!\!$ and~denoting by $\hat K^{\delta^{^G}}$ the converse image in 
$\hat C_G (Q)$ of the fixed~points of~$E_H(Q_{\delta^{^G}})$ in $\hat Z$ and by $K^{\delta^{^G}}$
the {\it $k^*\-$quotient\/} of~$\hat K^{\delta^{^G}}\!\!\,,$ as above $H\.K^{\delta^{^G}}\!\!$ is a normal subgroup of~$G$ and therefore, arguing by induction, we may assume that either $C_H (Q) =  K^{\delta^{^G}}\!\!$
or~$G= H\.K^{\delta^{^G}}\!\!\,.$

\smallskip 
In the first case, it follows from Proposition~£3.8 that 
$$F =  E_H (Q_{\delta^{^G}})\i E_G (Q_{\delta^{^G}}) = E
\eqno £3.13.61\phantom{.}$$
 so that $S$ is indeed $F\-$stable (cf.~£3.9). In the second case, since  we have (cf.~Proposition~£3.8)
 $$F/E_H (Q_{\delta^{^G}})\cong K^{\delta^{^G}}\!\!/K^{\delta}\cong G/H \cong Z
  \eqno £3.13.62,$$
 the number of points of $\O_*\hat Z$ coincides with its {\rm $\O\-$rank\/} which forces the $k^*\-$group  isomorphism $\hat Z\cong k^*\times Z\,;$ in particular, isomorphism~£3.13.54 becomes the $Q\-$interior algebra isomorphism
 $$(\O G)_\delta \cong (\O H)_\delta\, Z = \bigoplus_{z\in Z} (\O H)_\delta\.z
 \eqno £3.13.63\phantom{.}$$
 and therefore we have $(\O G)_\delta^Q \cong (\O H)_\delta^Q\, Z \,.$

 \smallskip
 Thus, since $Q = P\,,$ we may assume that the image $i$ of ${1\over \vert Z\vert}\;\sum_{z\in Z} z$ in 
 $(\O G)_\delta\i \O G$ belongs to $\delta^{^G} = \gamma^{^G}\!\!$ and then we get (cf.~£3.9.1)
 $$S\otimes_\O \O_*\hat L\cong i(\O G)i\cong (\O H)_\delta^Z
 \eqno £3.13.64.$$
 But, it follows from [10, Proposition~7.4] that there is a unique $j + J\big((\O H)_\delta^Q\big)\-$ conjugacy class of $k^*\-$group homomorphisms
 $$\hat\alpha : Q\rtimes\hat F\too \big((\O H)_\delta\big)^*
 \eqno £3.13.65\phantom{.} $$
 mapping $u\in Q$ on $u\.j\,;$ then, since $Z$ is a $p'\-$group, it follows from [3, Lemma~3.3 and~Proposition~3.5]
 that we can choose $\alpha$ in such a way that $Z$ normalizes $\alpha (\hat F)$ and then we have 
 $[Z,\alpha (\hat F)]\i k^*\,.$ In this case, $\alpha(\hat F)$ stabilizes $(\O H)_\delta^Z$  and, as a matter of fact, this $\O\-$algebra contains  $\alpha (\hat F)^Z = \hat E\,.$
 \eject

 \medskip
 Consequently, throughout isomorphisms~£3.13.64, $F$ acts on the $Q\-$interior algebra~$S\otimes_\O \O_*\hat L$
and therefore it acts on the quotient
 $$S\otimes_\O \O_*\hat L\big/J(S\otimes_\O \O_*\hat L)\cong S\otimes_\O k_*\hat E
 \eqno £3.13.66;$$
 but, $k_*\hat E$ has a trivial  $Q\-$interior algebra structure and it is semisimple, so that we have
 $$ k_*\hat E\cong \prod_\theta  (k_*\hat E)(\theta)
 \eqno £3.13.67,$$  
 where $\theta$ runs over the set of points of $k_*\hat E\,,$ and then any tensor product $T_\theta = S\otimes_\O (k_*\hat E)(\theta)$
 is a Dade $Q\-$algebra over $k$  {\it similar\/} to $S\otimes_\O k$ [11,~1.5];
hence, $F$ permutes these Dade $Q\-$algebras which amounts to saying that it fixes the {\it similarity class\/}  
of $S\otimes_\O k\,;$ finally, it follows from 
 [11, 1.5.2] that $S$ is also $F\-$stable. We are done.

 \bigskip
\bigskip
\noindent
{\bf £4\phantom{.} Normal sub-blocks of nilpotent blocks}
\bigskip

£4.1\phantom{.} With the notation of section 3, assume now that the block $b$ of $G$ is nilpotent; since
we already know that $(\O G)_\gamma\cong S\otimes_O \O P$ for a suitable Dade $P\-$algebra $S$
[9, Main Theorem], the block $b$ is also inertial and therefore we already have proved that the normal 
sub-block $c$ of $H$ is inertial too; let us show with the following example --- as a matter of fact, the example which
 has motivated this note --- that the block $c$ need not be nilpotent.

 \medskip
 \noindent
 {\bf Example~£4.2}\phantom{.} Let $\frak F$ be a finite field of characteristic different from $p\,,$ 
 $q$ the cardinal of $\frak F$ and $\frak E$ a field extension of $\frak F$ of degre $n\not= 1\,;$ denoting by $\Phi_n$
 the $n\-$th {\it cyclotomic polynomial\/}, assume that $p$ divides $\Phi_n (q)$ but not $q -1\,,$
that $\Phi_n (q)$ and $q -1$ have a nontrivial common divisor $r$ --- which has to be a prime number{\footnote{\dag}{\cds We thank Marc Cabanes for this remark.}} --- and that $n$ is a power of $r\,.$ For instance, the triple $(p,q,n)$
could be $(3,5,2)\,,$ $(5,3,4)\,,$ $(7,4,3)$ $\dots$
\smallskip
Set $G = GL_\frak F(\frak E)$ and $H = SL_\frak F(\frak E)\,,$ and respectively denote by $T$ and by $W$ 
the images in $G$ of the multiplicative group of $\frak E$ and of the Galois group of the extension $\frak E/\frak F\,;$
since $p$ does not divide $q-1\,,$ $T\cap H$ contains the Sylow $p\-$subgroup $P$ of $T$ and, since $p$ divides
$\Phi_n (q)\,,$ we have
$$C_G (P) = T\qq N_G (P) = T\rtimes W
\eqno £4.2.1;$$
consequently, since $W$ acts regularly on the set of generators of a Sylow $r\-$subgroup of $T\,,$ a generator $\varphi$ of the Sylow $r\-$subgroup of  ${\rm Hom}(T,\Bbb C^*)$ determines a local point $\gamma$ of $P$ on $\O G$ such that
$$N_G (P_\gamma) = T = C_G (P)
\eqno £4.2.2\phantom{.}$$
and, by the {\it Brauer First Main Theorem\/}, $P_\gamma$ is a defect pointed group of a block $b$ of $G$
which, according to [13, Proposition~5.2], is {\it nilpotent\/} by equality~£4.2.2.

\smallskip
On the other hand, since $r$ divides $q -1\,,$ the restriction $\psi$ of $\varphi$  to the intersection $T\cap H
= C_H (P)$ has an order
strictly smaller than $\varphi$ and therefore, since we clearly have 
$$N_H (P)/C_H (P)\cong W
\eqno £4.2.3,$$
$r$ divides $\vert N_H (P_\delta)/C_H (P)\vert$ where $\delta$ denotes the local point of $P$ on $\O H$ 
determined by $\psi\,;$ once again  by the {\it Brauer First Main Theorem\/}, $P_\delta$ is a defect pointed group of a block $c$ of $H\,,$ which is clearly a normal sub-block of the block $b$ of $G$ and it is {\it not\/} nilpotent
since $r$ divides $\vert N_H (P_\delta)/C_H (P)\vert\,.$

\bigskip
\noindent
{\bf Corollary~£4.3}\phantom{.} {\it A block $c$ of a finite group $H$ is a normal sub-block of
a nilpotent block $b$ of a finite group $G$ only if it is inertial and has an Abelian inertial quotient.\/}
\medskip
\noindent
{\bf Proof:} We already have proved that $c$ has to be inertial. For the second statement, we borrow the notation of Proposition~£3.8; on the one hand, since the block $b$ is nilpotent,  we know that $E_G (Q_{\delta^{^G}})$ is a 
$p\-$group; on the other hand, it follows from this proposition that $E_H (Q_{\delta^{^G}})$ is a normal subgroup
of $F$ and that $F/E_H (Q_{\delta^{^G}})$ is  Abelian; since the inertial quotient $F$ is a $p'\-$group, we have
$E_H (Q_{\delta^{^G}}) = \{1\}$ and $F$ is Abelian. We are done.

\medskip
\noindent
{\bf Remark~£4.4}\phantom{.} Conversely, if $P$ is a finite $p\-$group and $E$ a finite Abelian $p'\-$group acting faithfully on $P\,,$ the unique block of $\hat L = P\rtimes \hat E$ for any central $k^*\-$extension of $E$ is a normal sub-block of a nilpotent block of a finite group obtained as follows. Setting
$$Z = {\rm Hom}(E,k^*)
\eqno £4.4.1,$$
it is clear that $Z$ acts faithfully on $\hat E$ fixing the $k^*\-$quotient $E\,;$ thus, the semidirect product
$\hat E\rtimes Z$ still acts on $P$ and we finally consider the semidirect product
$$\hat M = P\rtimes (\hat E\rtimes Z) = \hat L\rtimes Z
\eqno £4.4.2.$$
Then, we clearly have
$$(\O_*\hat M)(P) \cong k\big(Z(P)\times Z\big)
\eqno £4.4.3\phantom{.}$$
and therefore any group homomorphism $\varepsilon\,\colon Z\to k^*$ determines a local point of $P$ on 
$\O_*\hat M$ --- still noted  $\varepsilon\,;$ but $E$ acts on $kZ\,,$ regularly permuting the set of its points;
hence, we get
$$N_{\hat M}(P_\varepsilon) = k^*\times P\times Z
\eqno £4.4.4\phantom{.}$$
and therefore $P_\varepsilon$ is a defect pointed group of the nilpotent block $\{1_{_{\O_*\hat M}}\}$ of~$\hat M\,.$

\vfill
\eject
\bigskip
\noindent
{\bf References}
\bigskip
\noindent
[1]\phantom{.} Michel Brou\'e and Llu\'\i s Puig, {\it Characters and Local
Structure in $G\-$alge-bras,\/} Journal of Algebra, 63(1980), 306-317.
\smallskip\noindent
[2]\phantom{.} Michel Brou\'e and Llu\'\i s Puig, {\it A Frobenius theorem for
blocks,\/} Inventiones math., 56(1980), 117-128.
\smallskip\noindent
[3]\phantom{.} Yun Fan and Lluis Puig, {\it On blocks with nilpotent
coefficient extensions\/}, Algebras and Representation Theory, 1(1998),
27-73 and Publisher revised form, 2(1999), 209.
\smallskip\noindent
[4]\phantom{.} Daniel Gorenstein, {\it ``Finite groups''\/} Harper's Series,
1968, Harper and Row.
\smallskip\noindent
[5]\phantom{.} James Green, {\it Some remarks on defect groups\/}, Math. Zeit.,
107(1968), 133-150.
\smallskip\noindent
[6]\phantom{.} Burkhard K\"ulshammer and Llu\'\i s Puig, {\it Extensions of
nilpotent blocks}, Inventiones math., 102(1990), 17-71.
\smallskip\noindent
[7]\phantom{.} Llu\'\i s Puig, {\it Pointed groups and  construction of
characters}, Math. Zeit. 176(1981), 265-292. 
\smallskip\noindent
[8]\phantom{.} Llu\'\i s Puig, {\it Local fusions in block source algebras\/},
Journal of Algebra, 104(1986), 358-369. 
\smallskip\noindent
[9]\phantom{.} Llu\'\i s Puig, {\it Nilpotent blocks and their source
algebras}, Inventiones math., 93(1988), 77-116.
\smallskip\noindent
[10]\phantom{.} Llu\'\i s Puig, {\it Pointed groups and  construction of
modules}, Journal of Algebra, 116(1988), 7-129.
\smallskip\noindent
[11]\phantom{.} Llu\'\i s Puig, {\it Affirmative answer to a question of Feit},
Journal of Algebra, 131(1990), 513-526.
\smallskip\noindent
[12]\phantom{.} Llu\'\i s Puig, {\it Alg\`ebres de source de certains blocks des groupes
de Chevalley\/}, in {\it ``Repr\'esentations lin\'eaires des groupes finis''\/}, Ast\'erisque, 181-182 (1990),
Soc. Math. de France
\smallskip\noindent
[13]\phantom{.} Llu\'\i s Puig, {\it Une correspondance de modules pour  les
blocks \`a groupes de d\'efaut ab\'eliens}, Geometri\ae\  Dedicata, 37(1991),
9-43. 
\smallskip\noindent
[14]\phantom{.} Llu\'\i s Puig, {\it On Joanna Scopes' Criterion of equivalence for blocks 
of symmetric groups}, Algebra Colloq., 1(1994), 25-55.
\smallskip\noindent
[15]\phantom{.} Llu\'\i s Puig, {\it ``On the Morita and Rickard
equivalences between Brauer blocks''\/}, Progress in Math., 178(1999), Birkh\"auser,
Basel.
\smallskip\noindent
[16]\phantom{.} Llu\'\i s Puig, {\it Source algebras of $p\-$central group 
extensions\/}, Journal of Algebra, 235(2001), 359-398.
\smallskip\noindent
[17]\phantom{.} Llu\'\i s Puig, {\it ``Blocks of Finite Groups''\/},
Springer Monographs in Mathematics, 2002, Springer-Verlag, Berlin, Barcelona.
\smallskip\noindent
[18]\phantom{.} Llu\'\i s Puig, {\it ``Frobenius categories versus Brauer blocks''\/}, Progress in Math., 
274(2009), Birkh\"auser, Basel.
\smallskip\noindent
[19]\phantom{.} Llu\'\i s Puig, {\it Block Source Algebras  in p-Solvable
Groups},  Michigan Math. J. 58(2009), 323-328

\end